\documentclass[12pt]{article}
\usepackage{amssymb,latexsym,amsmath}
\usepackage{graphicx}
\def\C{\mathbb C}
\def\R{\mathbb R}
\def\N{\mathbb N}

\newtheorem{thm}{Theorem}[section]
\newtheorem{lem}{Lemma}[section]

\newtheorem{prop}{Proposition}[section]

\begin{document}
\sffamily
\title{Bank-Laine functions, the Liouville transformation\\
and 
the Eremenko-Lyubich class}
\author{J.K. Langley}
\maketitle

\centerline{\textit{For Larry Zalcman, on the occasion of his retirement from Journal d'Analyse Math\'ematique}}

\begin{abstract}
The Bank-Laine conjecture concerning the oscillation of solutions of second order homogeneous linear 
differential equations has recently been disproved by Bergweiler and Eremenko.  It is shown here, however, that the conjecture is true
if the set of finite critical and asymptotic values of the coefficient function is bounded. 
It is also shown that if $E$ is a Bank-Laine function of finite order with infinitely many zeros, all real and positive, then its zeros must have exponent of convergence
at least $3/2$, 
and an example is constructed via quasiconformal surgery to demonstrate that this result is sharp.
%that there exists a  real Bank-Laine function of finite order, all of whose zeros are real and positive.
%these having exponent of convergence $3/2$. 
\end{abstract}

\section{Introduction}

If $f$ is a non-constant entire function, let
$$
\rho (f) = \limsup_{r \to + \infty} \frac{ \log^+ T(r, f)}{\log r} \, , \quad 
\lambda (f) = \limsup_{r \to + \infty} \frac{ \log^+ N(r, 1/f)}{\log r}  \leq \rho (f),
$$
denote its order of growth and the exponent of convergence of its zeros \cite{Hay2}. 
In their landmark paper \cite{BIL1}, 
Bank and Laine proved the  following results on the oscillation of solutions of 
%the equation
\begin{equation}
 \label{de1}
y'' + A(z) y = 0.
\end{equation}
%Bank and Laine proved the  following results. 

\begin{thm}[\cite{BIL1}]
 \label{thmA}
Let $A$ be an entire function, let $f_1, f_2$ be linearly independent solutions 
of (\ref{de1}) and let $E=f_1f_2$, so that $\lambda (E) = \max \{ \lambda (f_1), \lambda (f_2) \}$.\\
(i) If $A$ is a polynomial of degree $n > 0$ then $\lambda(E) = (n+2)/2 $.\\
%at least one $f_j$ has $\lambda(f_j) = \rho(f_j) = (n+2)/2$.\\
(ii) If $\lambda(E) < \rho(A) < + \infty$ 
%and $\lambda(E) < + \infty$ 
then $\rho(A) \in \N = \{ 1, 2, \ldots \}$.\\
(iii) If $A$ is transcendental and $\rho(A) < 1/2$ then $\lambda(E) = + \infty$.
%, so that at least one $f_j$ has $\lambda(f_j) = + \infty$. 
\end{thm}
The case where $1/2 \leq \rho(A) < 1$ was considered by Rossi \cite{Ros} and Shen \cite{Shen}.

\begin{thm}[\cite{Ros,Shen}]
 \label{thmB}
Let $A$ be an entire function of order $ \rho(A)$ and let $E = f_1f_2$, where 
$f_1, f_2$ are linearly independent solutions of (\ref{de1}). If $\rho(A) = 1/2$ then  $\lambda (E) = + \infty$, while 
\begin{equation}
 \label{roineq}
\frac1{\rho(A)} + \frac1{\lambda (E) } \leq 2 \quad \text{ if $1/2 < \rho(A) < 1$.}
\end{equation}
In particular, if $1/2 \leq \rho(A) < 1$ then $\rho(E)  > 1$. 
\end{thm}
The methods of \cite{BIL1} 
%and many subsequent papers \cite{Lai1} 
focused on the product
$E = f_1 f_2$ of linearly independent solutions $f_j$ of (\ref{de1}), and in particular  on 
the equation 
\begin{equation}
 \label{bleq}
 4A = \left( \frac{E'}{E} \right)^2 - 2 \, \frac{E''}{E} - \frac{c^2}{E^2} , \quad c = W(f_1, f_2) ,
\end{equation}
linking $E$ and $A$, in which the Wronskian  $ W(f_1, f_2) = f_1f_2' - f_1'f_2  $ is constant by Abel's identity. 
The paper \cite{BIL1} inspired much subsequent activity concerning
the zeros of solutions of (\ref{de1}) and, more generally,  linear differential equations with entire coefficients
\cite{Lai1},
and gave rise to 
%what became known as 
the Bank-Laine conjecture -- \textit{let $A$ be a transcendental entire function of
finite order  $\rho(A)$
and let $f_1, f_2$ be linearly
independent solutions of $(\ref{de1})$: if $\lambda (f_1f_2)$ is finite then  
$\rho (A) \in \N $}. 
However,  two remarkable recent papers of Bergweiler
and Eremenko 
\cite{bebl1,bebl2} show via  quasiconformal constructions  not only that the Bank-Laine conjecture is  false,
but also that  the inequality (\ref{roineq}) is sharp.

When $A$ is a non-constant polynomial in (\ref{de1}), satisfying $A(z) = a_n z^n (1+o(1))$ as $z \to \infty$, 
there are $n+2$  critical rays given by $\arg z = \theta^*$, where $a_n e^{i (n+2) \theta^*} $ is real and positive, and 
the Liouville transformation
\begin{equation}
 \label{ltdef}
Y(Z) = A(z)^{1/4} y(z) , 
\quad Z = \int_{z_1}^z \, A(t)^{1/2} \, dt ,
%= \int_{v_1}^v e^{w/2} \phi'(w) \, dw  
\end{equation}
may be applied in sectors symmetric about these rays. 
This reduces (\ref{de1}) to a sine-type equation 
%$Y'' + (1+o(1))Y = 0$
$$ 
\frac{d^2Y}{dZ^2} + \left(1 + \frac{O(1)}{Z^2} \right) Y = 0, 
$$
for which solutions asymptotic to $e^{\pm iZ}$ on a sectorial region in the $Z$ plane are delivered by Hille's method \cite{Hil1,Hil2}.
On one side of the critical ray, one  of the corresponding solutions $A(z)^{-1/4} e^{\pm iZ} (1+o(1))$
of (\ref{de1}) is large while the other is small, and
these roles are reversed as the critical ray is crossed. 

In contrast, for transcendental entire  $A$, although a
local analogue of Hille's method was developed in \cite{La5}, applying on small neighbourhoods of maximum modulus points of $A$,
%arising from the Wiman-Valiron theory \cite{Hay5}, 
the analytic continuation and estimation of $Z$ in (\ref{ltdef}) present substantial difficulties. However, it turns out that for a certain class
of entire functions $A$ the 
transformation (\ref{ltdef}) may be adapted so as to be readily applicable on components where $|A(z)|$ is large.

%For a transcendental entire function $A$, the singular values of the inverse function $A^{-1}$ are the asymptotic values (that is, values approached by $A(z)$ along paths tending to infinity) and the critical values. 
The Eremenko-Lyubich class $\mathcal{B}$ plays a key role in complex dynamics  \cite{Ber4,EL,sixsmithEL} and
consists of those 
transcendental meromorphic functions $A$ with the following property:
there exists a positive real number $M = M(A)$ such that all finite critical 
and asymptotic values of $A$ have modulus less than $M$. Now suppose that $A \in \mathcal{B}$ is entire. Then, by standard results from \cite[p.287]{Nev}
(see also \cite{BE}),
all components $U_M$ of the set $\{ z \in \C : |A(z)| > M \}$ correspond to 
logarithmic singularities of $A^{-1}$ over $\infty$; in particular, $v = \log A(z)$ maps each such $U_M$ conformally onto the half-plane $H$ given by
${\rm Re} \, v > \log M$. Under the change of variables
\begin{equation}
 \label{b4}
 A(z) = e^v, \quad z = \phi (v) , \quad \frac{A'(z)}{A(z)} = \frac{dv}{dz} = \frac1{\phi'(v)} ,
\end{equation}
in which  $z = \phi(v)$ is the inverse mapping from $H$ to $U_M$, 
a solution $y(z)$ of (\ref{de1}) on $U_M$ transforms to a solution $w(v) = y(z)$ on $H$ of 
\begin{equation}
  \label{p2}
  w''(v) - \frac{\phi''(v)}{\phi'(v)} w'(v) + e^v \phi'(v)^2 w(v) = 0,
 \end{equation} 
and the second formula in (\ref{ltdef}) becomes, for a suitable choice of  $z_1 = \phi(v_1)$, 
\begin{equation}
 \label{ltvplane}
Z = \int_{v_1}^v e^{u/2} \phi'(u) \, du .
\end{equation}
%for $z \in U_M$, $v \in H$ and .
%, in which $\phi'(w)$  varies relatively slowly on $H$ by classical results on conformal mappings. 
The fact that $\phi'$ varies relatively slowly on $H$, by classical theorems on conformal mappings \cite{Hay9}, makes it possible
%se transformations will be used 
to prove the following theorem. 

\begin{thm}
 \label{thmel}
 Suppose that $A$ is a transcendental entire function in the Eremenko-Lyubich class $\mathcal{B}$, and let $E = f_1 f_2$,
where $f_1, f_2$ are linearly independent solutions of (\ref{de1}). Then exactly one of the following holds.
\\
(A) The functions  $A$ and  $E$ satisfy  $\rho(A) = \rho (E)  = 1$ and 
 \begin{equation}
  \label{AEest}
  T(r, A) + T(r, E) = O(r) \quad \hbox{as} \quad r \to + \infty. 
 \end{equation}
(B) There exists $d > 0$ such that the  zeros of $E$ satisfy 
\begin{equation}
 \label{Ezerosconc}
n(r, 1/E) > \exp \left( d r^{1/2} \right) \quad \hbox{as} \quad r \to + \infty ,
\end{equation}
and in particular  $\rho (E) = \lambda(E) =  + \infty$. 
\end{thm}
It follows from Theorem \ref{thmel} that  the Bank-Laine conjecture, 
despite being false in general \cite{bebl1}, is true when 
the coefficient $A$ is entire and in the class $\mathcal{B}$.
%The estimate (\ref{AEest}) is sharp, 
An example going back to \cite{BIL1} shows that each of conclusions (A) and (B) can occur: if $A(z) = - e^{2z} -1/4 $ then  (\ref{de1}) has solutions
$$
f_1(z) = e^{-z/2} \exp \left( - e^z \right), \quad 
f_2(z) = e^{-z/2} \exp \left( e^z \right), \quad  f_1(z) f_2 (z) =  e^{-z} , \quad \rho(f_1f_2) = 1, 
$$
as well as solutions 
$$
g_1(z) = e^{-z/2} \sinh \left(  e^z \right), \quad 
g_2(z) = e^{-z/2} \cosh \left( e^z \right), \quad  \lambda(g_1 g_2) = + \infty .
$$
An example will be given in Section \ref{pfthme1} to show that the exponent $1/2$ in (\ref{Ezerosconc}) is sharp. 
%Then the $f_j$ solve (\ref{de1}), in which $A(z) = - e^{2z} -1/4 $ evidently belongs to $\mathcal{B}$, while $\rho(A) = \rho (E)  = 1$ and $\lambda(E)=0$. 

The second main result of this paper concerns the location of zeros of Bank-Laine functions, that is, entire functions $E$ such
that $E(z)=0$ implies $E'(z) = \pm 1$. By \cite[Lemma C]{BIL3},
an entire function $E$ is a Bank-Laine function if and only if 
$E=f_1f_2$, where $f_1, f_2$ are linearly independent solutions of (\ref{de1}) with $A$ entire and $W(f_1, f_2) = 1$. 
Although 
a Bank-Laine function with no restriction on its growth
may have an arbitrary sequence $(a_n)$ of zeros, subject only to $a_n \to \infty$ without repetition \cite{Shen2},
the following result was proved in \cite{qcsurg} concerning Bank-Laine functions with real zeros.

\begin{thm}[\cite{qcsurg}]
 \label{thmDL}
Let $E$ be a Bank-Laine function of finite order, with infinitely many zeros, all real, and denote by $n(r)$ the number of zeros of $E$
lying in $[-r, r]$.  Then $n(r) \neq o(r)$ as $r \to + \infty$.  If, in addition,  all zeros of $E$ are positive, then 
$n(r) \neq O(r)$ as $r \to + \infty$. 
\end{thm}

The first assertion of Theorem \ref{thmDL} is evidently sharp, because of $\sin z$. The next theorem will establish a sharp lower bound for $\lambda(E)$ 
when $E$ is a Bank-Laine function of finite order with infinitely many zeros, all real and positive. Here it is sufficient to consider the case where $E$ is real entire, because 
otherwise it is possible to write $E = \Pi  e^{P+iQ}$, where $\Pi$ is the canonical product over the zeros of $E$,
while  $P$ and $Q$ are real polynomials; thus $e^{iQ(z)} = \pm 1$ at every zero of $E$ 
and $F = \Pi e^P$ is also a Bank-Laine function.

\begin{thm}
 \label{thmBLpositive}
Let $E$ be a real Bank-Laine function of finite order, with infinitely many zeros, all real and positive. Then the exponent of convergence $\lambda(E)$ of the zeros of
$E$ is at least $3/2$. 
%If, in addition, $E$ is real entire and
Moreover, if 
$\lambda(E) = 3/2$ then $E$ and the associated coefficient function $A$ have order $\rho(E) = \rho (A) = 3/2$. 
\end{thm}

To demonstrate the sharpness of Theorem \ref{thmBLpositive},
quasiconformal techniques  will be used  in Section
%s \ref{example1},  \ref{example2} and  
\ref{example3}
to  construct a real Bank-Laine function $E$, with only positive zeros, such that 
$E$ and its associated coefficient function $A$ satisfy 
$\lambda(E) = \rho(E) = \rho(A) = 3/2$, so that $A$
%these having exponent of convergence $3/2$. The coefficient function $A$  associated with $E$ also has order $3/2$, and so 
provides a further counter-example to the Bank-Laine conjecture. 
%shown further 

The author thanks the referee for an extremely careful reading of the manuscript and for numerous helpful suggestions. 

\section{A refinement of Hille's method }\label{hille}

The following lemma is an extension of a method from \cite{La5}, and provides  bounds for the error terms in Hille's method \cite{Hil1,Hil2}. 

\begin{lem}
 \label{lemhil1}
Let $c > 0$ and $0 < \varepsilon < \pi $. Then there exists
$d > 0$, depending only on $c$ and $\varepsilon$, with
the following properties. 
Suppose that the function $A$ is analytic, with $|1-A(z)| \leq c |z|^{-2}$,
on a domain containing 
\begin{equation*}
\Omega = \Omega_{R,S} =
\{ z \in \C : \, 1 \leq R \leq |z| \leq S < + \infty , \, | \arg z | \leq \pi - \varepsilon \}.
\label{H2}
\end{equation*}
%Set $A(z) = 1-F(z)$. 
Then the equation (\ref{de1})
%\begin{equation}w'' + ( 1 - F(z) ) w = 0\label{H3}\end{equation}
has linearly independent solutions $U(z), V(z)$ satisfying
\begin{eqnarray}
U(z) &=& e^{-iz} ( 1 + \delta_1(z) ), \quad
U'(z) = -i e^{-iz} ( 1 + \delta_2(z) ), \quad \nonumber \\
V(z) &=& e^{iz} ( 1 + \delta_3(z) ), \quad
V'(z) = i e^{iz} ( 1 + \delta_4(z) ), \quad
\label{H4}
\end{eqnarray}
in which 
\begin{equation}
| \delta_j (z) | \leq \frac{d}{ |z|} \quad \hbox{for} \quad z \in \Omega_{R,S}^* = \Omega_{R,S} \setminus
\{ z \in \C : \, {\rm Re } (z) < 0 , \, | {\rm Im }  (z) | < R \} .
\label{H5}
\end{equation}
\end{lem}
\textit{Proof.}
%Here $\Omega_{R,S}^*$ can be thought of as $\Omega$ with the ``shadow'' of $D(0, R)$ removed.
Let $X = S e^{i \sigma } $, where
$\sigma = \min \{ \pi /2 , \pi - \varepsilon \}$. Choose an analytic solution
$v$ on $\Omega$ of 
\begin{equation}
v'' + 2i v' - F v = 0, \quad F = 1-A, 
\label{H6}
\end{equation}
such that $v(X) = 1, v'(X) = 0$, and write
\begin{equation}
L(z) = v(z)  - 1 + \frac{1}{2i} \int_X^z ( e^{2i(t-z)} - 1 ) F(t) v(t) \, dt ,
\quad L'(z) = v'(z) - \int_X^z e^{2i(t-z)} F(t) v(t) \, dt ,
\label{H7}
\end{equation}
so that 
$$
L''(z) = v''(z) + 2i \int_X^z e^{2i(t-z)} F(t) v(t) \, dt
- F(z) v(z) = - 2i L'(z) .
$$
Since $L(X) = L'(X) = 0$, the existence-uniqueness theorem
implies that $L(z) \equiv 0$ on $\Omega$. 

Now let $z \in \Omega_{R,S}^*$ and let $\gamma_z$  describe the clockwise arc of the circle $|t| = S$ from $X$  to the first point $x$
of intersection with the line ${\rm Im} (t) = {\rm Im }  (z)$,
followed by the straight line segment
from $x$ to $z$; then 
%$ {\rm Im }  (t- z) \geq 0$ and hence
$| e^{ 2i (t-z) } | \leq 1$ on $\gamma_z \subseteq \Omega$. 
Since $L(z) = 0$,
(\ref{H7}) gives
\begin{equation}
| v(z) - 1 | \leq \int_X^z | F(t) v(t) | \, |dt|, \quad
| v(z) | \leq 1 + \int_X^z | F(t) v(t) | \, |dt| . 
\label{H9}
\end{equation}

Now parametrize $\gamma_z$ by $t = \zeta(s)$, where 
$s$ denotes arc length on $\gamma_z$. Using (\ref{H9}), write 
$$
H( s ) =  
% 1 + \int_X^{\zeta(s)} | F(t)v(t) | \, |dt| =
1 + \int_0^s | F( \zeta (\sigma ) ) v ( \zeta (\sigma) ) | \, d\sigma ,
\quad 
%\frac{dH}{ds} 
H'(s) =  | F( \zeta(s))v( \zeta(s)) | \leq |F(\zeta(s))| H(s) ,
%\zeta \in \gamma_z .
$$
and
%using the second estimate of .Thus (\ref{H9}) delivers, as in the standard proof of Gronwall's lemma \cite{Hil2}, 
$$
| v(\zeta(s)) - 1 | \leq H(s) - 1 =  \exp \left( \int_0^s \frac{H'(\sigma)}{H(\sigma)}  \, d\sigma \right) - 1 \leq
 \exp \left( \int_0^s | F( \zeta (\sigma) )  | \, d\sigma \right) - 1 ,
$$
%\begin{eqnarray*} | v(\zeta(s)) - 1 | &\leq& H(s) - 1 = \frac{H(s)}{H(0)} - 1   \\ &\leq&  \exp \left( \int_0^s \frac{H'(\sigma)}{H(\sigma)}  \, d\sigma \right) - 1 \\
%&\leq& \exp \left( \int_0^s | F( \zeta (\sigma) )  | \, d\sigma \right) - 1 , \end{eqnarray*}
%$$| v(\zeta(s)) - 1 | \leq H(s) - 1 = \frac{H(s)}{H(0)} - 1   \leq \exp \left( \int_0^s | F( \zeta (\sigma) )  | \, d\sigma \right) - 1 ,$$
which leads to
\begin{equation}
 \label{H10}
|v(z)-1| \leq \exp \left( I_z   \right) - 1, \quad I_z = \int_X^z | F(t) | \, |dt|  .
\end{equation}
%\begin{eqnarray}|v(z)-1| &\leq& \exp \left( I_z   \right) - 1, \nonumber \\%&=& \exp \left( I_z   \right) - 1, \nonumber \\
%I_z &=&  \int_0^{|\gamma_z|} | F( \zeta (\sigma) )  | \, d\sigma = \int_X^z | F(t) | \, |dt|  .\label{H10}\end{eqnarray}
Let $d_1 , d_2 , \ldots $ denote positive constants which depend
only on $c$
and $\varepsilon$. The circle $|t| = S$ 
contributes at most $d_1 S^{-1} \leq d_1 |z|^{-1}$ to $I_z$
in (\ref{H10}), while  the contribution $J_z$ from the horizontal part of $\gamma_z$ satisfies:
%is, if $| \arg z | \leq \pi /4$, at most
\begin{eqnarray*}
 J_z &\leq&  \int_{{\rm Re} \, z}^{+ \infty } \frac{c }{t^{2}} \,  dt = \frac{c}{ {\rm Re } \, z }  \leq \frac{d_2}{|z| } \quad \text{if $| \arg z | \leq \pi /4$;}\\
J_z &\leq& \int_{\R} \frac{c}{x^2 + ({\rm Im}\, z)^2 } \, dx \leq \frac{d_3}{ |{\rm Im} \, z |} \leq \frac{d_4}{ |z|} \quad \text{if $ \pi /4  \leq | \arg z | \leq 
\pi - \varepsilon $.}
\end{eqnarray*}
Since $R \geq 1$, (\ref{H7}) and (\ref{H10}) now deliver
$$
| v(z) - 1 | \leq \exp \left( \frac{d_5}{|z|} \right) - 1 \leq \frac{d_6}{ |z|} \leq d_6 ,
%$$using the fact that $R \geq 1$, and (\ref{H7}) gives$$
\quad 
| v'(z) | \leq \int_X^z | F(t) | (1+ d_6) \, |dt|  \leq \frac{d_7}{ |z|} .
$$
Now  set $V(z) = v(z) e^{iz} $; then (\ref{H6}) implies
that $V$ solves (\ref{de1}),   and the estimates (\ref{H4}) and (\ref{H5}) for $V$ follow at once.
To obtain $U$ it is only necessary to apply the above argument to the equation solved by $\overline{y(\bar z)}$ for every solution $y(z)$
of (\ref{de1}).
\hfill$\Box$
\vspace{.1in}

Unbounded sectorial regions may be handled as follows. 

\begin{lem}
 \label{lemhil2}
Suppose that $c > 0$ and $0 < \varepsilon < \pi $, and that the function $A$ is analytic, with $|1-A(z)| \leq c |z|^{-2}$, on  $\Omega' =
\{ z \in \C : 1 \leq R \leq |z| < + \infty , | \arg z | \leq \pi - \varepsilon \}$. 
%Set $A(z)=1-F(z)$. 
Then there exist $d > 0$, depending only on  $c$ and $\varepsilon$, and solutions $U, V$ of (\ref{de1})
on 
$$
\Omega'' =
\{ z \in \C :  R < |z| < + \infty , | \arg z | <  \pi - \varepsilon \}
\setminus \{ z : {\rm {\rm Re } } (z) \leq 0, |{\rm Im }  (z) | \leq R \} ,
$$
such that $U$ and $V$  satisfy $W(U, V) = 2i$ and (\ref{H4}), with $|\delta_j(z)| \leq d/|z|$, on $\Omega''$. 
\end{lem}
\textit{Proof.}
Taking a sequence $S_n \to +\infty $ yields solutions $U_n, V_n$ of (\ref{de1}) on $\Omega_{R,S_n}^*$, with  corresponding error terms  
$\delta_{j,n} (z), j = 1, 2, 3, 4$.
%$$\{ z \in \C : R \leq |z| \leq S_n, | \arg z | \leq \pi - \varepsilon \}\setminus \{ z : {\rm Re } (z) < 0, | {\rm Im }  (z) | < R \} ,$$
Here the functions $z \delta_{j,n} (z) $ are uniformly bounded, since the constant
$d$ is independent of $S$ in (\ref{H5}). Thus,
by normal families, it may be assumed  that the $U_n, V_n , \delta_{j, n}$ converge locally uniformly on $\Omega''$. 
The limit functions $U, V$ satisfy (\ref{H4}), with $|\delta_j(z)| \leq d/|z|$ on $\Omega''$. Since $W(U, V)$ is constant, by Abel's identity, 
it follows that $W(U, V) = 2i$.  
\hfill$\Box$
\vspace{.1in}

Finally, a change of variables $z \to -z$ shows that Lemmas \ref{lemhil1} and \ref{lemhil2}
hold if $\Omega_{R,S} $ and $ \Omega_{R,S}^* $, and correspondingly $\Omega'$ and $ \Omega''$,
are replaced by their reflections across the imaginary axis.

\section{Estimates in a half-plane}
Throughout this section let $H = \{ v \in \C : \text{Re} \, v >  0 \}$ and let $\phi : H \to \C \setminus \{ 0 \}$ 
be analytic and univalent.
%conformal mapping. 
For $v, v_1 \in H$, define $Z = Z(v, v_1)$ as in (\ref{ltvplane}) by
%with ${\rm Re} \, v_0 = 1$ define $Z$ by 
\begin{equation}
 \label{h1}
 Z(v, v_1) = \int_{v_1}^v e^{u/2} \phi'(u) \, du = 
2 e^{v/2} \phi'(v) - 2 e^{v_1/2} \phi'(v_1) - 2 \int_{v_1}^v e^{u/2} \phi''(u) \, du .
\end{equation}
Since $0 \not \in \phi (H)$
the image of $H$ under $\log \phi $ contains no disc of radius greater than $\pi$; thus 
applying Bieberbach's theorem and Koebe's one quarter theorem 
\cite[Theorems 1.1 and 1.2]{Hay9}  to $ \phi$ and $\log \phi$ respectively gives, for $u \in H$, 
\begin{equation}
 \label{h3}
\left| \frac{\phi''(u)}{\phi'(u)} \right| \leq \frac4{{\rm Re} \, u}, \quad  \left| \frac{\phi'(u)}{\phi(u)} \right| \leq \frac{4 \pi}{{\rm Re} \, u}  .
\end{equation}
The fact that the estimates (\ref{h3}) are independent of $\phi$ is the key to the results of this section and the proof of Theorem \ref{thmel}. 

\begin{lem}
 \label{lemfirstest}
Let $\varepsilon $  be a small positive real number. Then there exists a large positive real number $N_0$, depending  on $\varepsilon$ but not on $\phi$,
with the following property.

Let $v_0 \in H$ be such that  $S_0 = {\rm Re} \, v_0 \geq N_0$,
and define $v_1, v_2, v_3, K_2$ and $ K_3$ by 
\begin{equation}
\label{vjdef}
v_j = \frac{2^j S_0}{128} + i T_0, 
% v_j = v_0 - \frac{ {\rm Re} \, v_0 }{2^j} ,  
 \quad T_0 = {\rm Im} \, v_0, \quad 
K_j = \left\{ v_j + r e^{i \theta} : \, r \geq 0, \, - \frac{\pi}{2^j} \leq \theta \leq  \frac{\pi}{2^j} \right\}.
%, \quad j = 1, 2, 3.
%\ldots .
\end{equation}
Then the following three conclusions all hold:\\
(i)  $Z = Z(v, v_1)$ satisfies,  for $v \in K_2$,
%re exists $N > 0$ such that if ${\rm Re} \, v_0 > N$ then 
\begin{equation}
 \label{h2}
 Z =
Z(v, v_1 ) = \int_{v_1}^v e^{u/2} \phi'(u) \, du  = 2 e^{v/2} \phi'(v) (1 + \delta (v) ), \quad | \delta (v) | < \varepsilon .
\end{equation}
(ii) $\psi = \psi (v, v_1) = \log Z(v, v_1) $
is univalent on a domain containing $K_3$.\\
(iii) There exists a domain $D$, with $v_0 \in D \subseteq K_3$, mapped univalently by $Z$ onto 
a sectorial region $M_3$ satisfying
%which contains $ Z_0 = Z(v_0, v_1)  $ and is given by
\begin{equation}
 \label{Omegaimage}
Z_0 = Z(v_0, v_1) \in  M_3 = 
\{ Z \in \C : \, |Z_0|/8 < |Z| < + \infty, \, | \arg (\eta Z) |< 3 \pi /4 \} ,
%, \quad \log R_0 \geq \frac12 ( S_0 + \log  | \phi'(v_0) |   ) ,
\end{equation}
where $\eta = 1$ if ${\rm Re} \, Z_0 \geq 0$ and $\eta = -1$ if ${\rm Re} \, Z_0 < 0$.
\end{lem}
\textit{Proof.} To prove (i) assume that $S_0 = {\rm Re} \, v_0$ is large and let $v \in K_2$, so that
\begin{equation}
 \label{Tdef}
S = {\rm Re} \, v \geq \frac{S_0}{32} = 2 \, {\rm Re} \, v_1 .
\end{equation}
Now
$v_1$ may be joined 
to $v$ by a straight line segment $L_v$ which is parametrised with respect to $s = {\rm Re} \, u$, and an 
elementary arc length estimate $ |du| \leq ( \sec \pi/4) ds \leq 2 \, ds$ holds on $L_v$. 
Thus (\ref{h3}) delivers, for $u \in L_v$,  
\begin{equation}
 \label{h4}
 % \left| \log \frac{\phi'(u)}{\phi'(v) \right| \leq 
 | \phi'(u)| \leq \left( \frac{S}{s} \right)^8 | \phi'(v)|, \quad
 | \phi''(u)|\leq \frac4{s} \,| \phi'(u)| \leq \frac4{s} \, \left( \frac{S}s \right)^8 | \phi'(v)| ,
\end{equation}
which implies by (\ref{Tdef}) that 
\begin{equation}
 \label{h5}
\left| \frac{ e^{v_1/2} \phi'(v_1) }{ e^{v/2} \phi'(v) } \right| \leq
\left( \frac{S}{{\rm Re} \, v_1 } \right)^8 \exp \left( \frac12  {\rm Re} \, (v_1-v) \right)  
\leq S^8 \exp \left( -S/4 \right)  < \frac{\varepsilon}4 
\end{equation}
provided $ S_0$ is large enough. Moreover, (\ref{h4})  leads to 
\begin{equation}
 \label{h6}
\left| \frac{ 1 }{ e^{v/2} \phi'(v) } \int_{v_1}^v e^{u/2} \phi''(u) \, du \right| \leq  \Psi(S) =
\frac{8 S^8}{  e^{S/2}} \int_1^S e^{s/2} s^{-9} \, ds   .
\end{equation}
Since $\lim_{S \to + \infty} \Psi (S) = 0$ by L'H\^opital's rule, (\ref{Tdef}) implies that
$\Psi(S) < \varepsilon/4$ if  
$  S_0$ is large enough.
Thus (\ref{h2}) follows from  (\ref{h1}), (\ref{h5}) and (\ref{h6}), which proves (i).

Next, (\ref{h2}) gives, on $K_2$,
$$
\psi (v)  = \psi (v, v_1) = \log Z(v, v_1) = \frac{v}2 + \log 2 + \log \phi'(v) + \delta_1(v), \quad | \delta_1(v) | \leq 2 | \delta (v) | < 2 \varepsilon .
$$
Since $\varepsilon$ is small and $S_0$ is large,
%Provided $S_0$ is large enough, 
(\ref{h3}), (\ref{vjdef})  and Cauchy's 
estimate for derivatives now deliver
\begin{equation}
 \label{psi'est}
 \left|  \psi'(v) - \frac12 \right| \leq \frac{8}{{\rm Re} \, v } \leq \frac14 ,
\end{equation}
and hence ${\rm Re} \, \psi'(v) > 0$, 
on a convex domain containing $K_3$, which proves (ii). 

Now let 
$$
L_3 = \{ v \in K_3 : \, {\rm Re} \, v \geq S_0/8 \} .
$$
%and, for eac $v \in L_3$, integrate 
Then,
%$L_0 \subseteq K_3$ and,
for $v \in L_3$, integration along the line segment from $v_0 $ to ${\rm Re} \, v + i T_0 $
followed by that from ${\rm Re} \, v + i T_0 $ to $v$ yields, in view of (\ref{psi'est}), 
\begin{equation}
 \label{psiest}
 \psi (v) - \psi (v_0) = \frac{v-v_0}2 + \eta(v), 
 \quad | \eta (v)| \leq 8 \left(  \left|  \log \frac{{\rm Re} \,v}{S_0} \right| +  \tan \frac\pi8  \right) .
%C_1 ,
\end{equation}
%where $C_1 > 0$ does not depend on $S_0$. 
Since $S_0$ is large this implies that,  for $v \in \partial L_3$ with ${\rm Re} \, v = S_0/8$,
$$
{\rm Re} \, (\psi (v) - \psi (v_0) ) \leq   - \, \frac{7 S_0}{16} + 8 \left( \log 8 + \tan \frac\pi8  \right)  \leq \log \frac1{16} .
$$
On the other hand,  all other $v \in \partial L_3$  satisfy, by (\ref{vjdef}) and (\ref{psiest}), 
\begin{eqnarray*}
 | {\rm Im} \, (v-v_0) | &\geq&  \left( {\rm Re} \, v - \frac{S_0}{16} \right)  \tan \frac\pi8 \geq \frac{ {\rm Re} \, v}2 \,  \tan \frac\pi8  \, , \\
| {\rm Im} \, (\psi (v) - \psi (v_0) )| &\geq&  \frac{ {\rm Re} \, v}4 \,  \tan \frac\pi8 - 8 \left( \left| \log \frac{{\rm Re} \,v}{S_0} \right| +  \tan \frac\pi8  \right) \geq  4 \pi .
\end{eqnarray*}
Moreover,    
${\rm Re} \, (\psi (v) - \psi (v_0) ) \to + \infty$ as $v \to \infty$ in $K_3$, again by (\ref{psiest}). 
Thus the strip 
$$
\left\{  \psi (v_0)  + \sigma+ i \tau  : \, \sigma \geq \log \frac18 \, , \, -  2 \pi \leq  \tau \leq 2 \pi \right\} 
$$
lies in the interior of $\psi(L_3)$, which completes the proof of (iii) and the lemma.
%, since the image of $L_3$ under the univalent function $ \psi = \log Z   $ must contain 
%Because $Z = \exp( \psi(v))$, this proves (iii). 
\hfill$\Box$
\vspace{.1in}

\begin{prop}
 \label{mainprop}
 There exists a positive  real number $N_1 $, independent of $\phi$,  with the following property. If $v_0 \in H$ satisfies
 \begin{equation*}
  \label{p1}
  \min \{  S_0 , \, |e^{v_0/2} \phi'(v_0) | \} > N_1 , \quad S_0 =  {\rm Re} \, v_0 ,
 \end{equation*}
 and if $w_1, w_2 $ are linearly independent solutions of (\ref{p2}) 
  with 
 \begin{equation}
  \label{p3}
  W(w_1, w_2) = \pm \phi' , \quad |w_1(v_0) w_2(v_0) | \geq 1,
 \end{equation}
 then $w_1w_2$ has a sequence of distinct zeros $\zeta_m \to \infty$ in $H$ which satisfy
%, for some $d_1 > 0$, 
 \begin{equation}
  \label{p4}
 | \phi( \zeta_m ) | = O  \left( \log m \right)^{2 }  \quad \hbox{as} \quad m \to + \infty .
\end{equation}
\end{prop}
\textit{Proof.} Observe first that, by Abel's identity, the Wronskian of any two local solutions of 
(\ref{p2}) is a constant multiple of $\phi'$. Fix a small positive $\varepsilon$ and assume that $v_0 \in H$,
that $w_1, w_2 $ are linearly independent solutions of
(\ref{p2}) which satisfy (\ref{p3}), and finally that  $S_0 $ and $ |e^{v_0/2} \phi'(v_0) |$ are both large. 
Let $v_1, v_2, v_3, K_2$ and $K_3$ be as in (\ref{vjdef}),  and define $Z$ 
%and $\psi = \log Z$ 
by (\ref{h1}). 
%Since $S_0$ is large, 
By Lemma \ref{lemfirstest}, $Z_0 = Z(v_0, v_1)$ is large and there exist $\eta \in \{ -1, 1 \}$ and
%, the choice of $\eta$ depending on ${\rm Re} \, Z_0$,  and
a domain $D \subseteq K_3$, both as in conclusion (iii), so that $M_3 = Z(D)$ satisfies (\ref{Omegaimage}). 
The change of variables 
\begin{equation}
 \label{p5}
 w(v) = e^{-v/4} W(Z), \quad w_j(v) = e^{-v/4} W_j(Z), 
%\quad v \in Z^{-1} (M_4)  ,
\end{equation}
transforms (\ref{p2}) on $D$ to the equation on $M_3$ given by
\begin{equation}
 \label{p6}
 W''(Z) + (1+G(Z)) W(Z) = 0, \quad G(Z) = \frac1{16 e^v \phi'(v)^2} \left( 1 + 4 \, \frac{\phi''(v)}{\phi'(v)} \right) .
\end{equation}
Here the derivatives in the first equation are with respect to $Z$, and
\begin{equation}
 \label{p7}
 |G(Z)| \leq \frac1{|Z|^2} 
\end{equation}
on $M_3 = Z(D)$, by  (\ref{h3}), (\ref{h2}) and the fact that  $S_0 = {\rm Re} \, v_0$ is large. Now apply Lemma \ref{lemhil2} 
with  
\begin{equation*}
 \label{Omeganewimage1}
 \Omega' = 
\{Z  \in \C : \, |Z_0|/4 \leq |Z| < + \infty, \, | \arg (\eta  Z) | \leq  5\pi /8 \} \subseteq M_3 , 
%, \quad \log R_0 \geq \frac12 ( S_0 + \log  | \phi'(v_0) |   ) ,
\end{equation*}
and let $M_4 = \Omega''$, so that  $Z_0 = Z(v_0, v_1) \in M_4 \subseteq \Omega' \subseteq M_3$, by the choice of $\eta$.
% in Lemma \ref{lemfirstest}.  
Since $|Z_0|$ is large, 
there exist solutions $U_1(Z), U_2(Z)$ of (\ref{p6}) on 
$M_4$, 
which satisfy $W(U_1, U_2) = 2i$ and 
\begin{equation}
\label{p8}
 |U_1(Z)e^{iZ} - 1| + | U_2(Z)e^{-iZ} - 1| \leq \frac{d}{|Z|} ,
\end{equation}
in which the  positive  constant $d$ is independent of  $v_0$ and $Z_0$, by (\ref{p7}).

Suppose first that, on $M_4$,
$$
W_1(Z) = \sigma_1 U_1(Z) , \quad W_2(Z) = \sigma_2 U_2(Z) , \quad \sigma_j \in \C \setminus \{ 0 \} .
$$
Then (\ref{h2}), (\ref{p3}) and (\ref{p5}) give
$$
\pm \phi' = W(w_1, w_2) = e^{-v/2} W(W_1, W_2) \frac{dZ}{dv} = W(W_1, W_2) \phi' = 2i \sigma_1 \sigma_2 \phi' ,
$$
so that $| \sigma_1 \sigma_2 | = 1/2$. But ${\rm Re} \, v_0$ and 
$|Z_0|$ are large, which implies, in view of (\ref{p5}) and (\ref{p8}),  that 
$$
w_1(v_0) w_2(v_0) = e^{-v_0/2} W_1(Z_0) W_2(Z_0) = e^{-v_0/2} \sigma_1 \sigma_2 U_1(Z_0) U_2(Z_0)
$$
is small,  a contradiction. 

Because $w_1$ and $w_2$ are interchangeable, it now follows  that at least one of $W_1$ and $W_2$, without loss
of generality $W_1$, is a non-trivial linear combination 
\begin{equation}
 \label{p9}
W_1(Z) = 
A_1 U_1(Z)- A_2 U_2(Z)  , \quad A_1, A_2 \in \C \setminus \{ 0 \},
\end{equation}
of $U_1, U_2$ on $M_4$. Fix a small positive $\kappa$ and suppose that 
$$
Z^*  = \frac1{2i} \log \frac{A_1}{A_2} +  \pi n , 
$$
where $n$ is an integer of large modulus and appropriate sign, depending on $\eta $.
Then $Z^* \in M_4$ and (\ref{p8})  implies that, on $|Z-Z^*| = \kappa$, 
$$
\frac1{2i} 
\log \frac{A_2 U_2(Z) }{A_1 U_1(Z)} - \pi n  = Z - Z^* + J(Z) , \quad |J(Z)| < \kappa .
$$ 
Hence $W_1$ has a zero  $Z^{**}$ with $|Z^{**}-Z^*| < \kappa $, by Rouch\'e's theorem and (\ref{p9}).

It follows that $W_1(Z)$ has distinct zeros $X_1, X_2, \ldots $,  which tend to infinity in $M_4$ and satisfy
$|X_m| \leq c_0 + c_1 m$, where $c_0, c_1, \ldots $ denote positive  constants which may depend on $v_0$ and $\phi$,
but not on $m$. By (\ref{h2}), these zeros $X_m$ satisfy, with $\zeta_m \in K_3$  and $\varepsilon$ small, 
$$
X_m = Z(\zeta_m) 
%= e^{\psi(\zeta_m) } 
= 2 e^{\zeta_m/2} \phi'(\zeta_m) (1 + \delta ( \zeta_m) ) , \quad | \delta ( \zeta_m) ) | < \varepsilon .
$$
Using (\ref{h3}) to estimate $| \log | \phi'(\zeta_m) | |$ then  gives, in view of (\ref{vjdef}), 
\begin{eqnarray*}
 |\zeta_m| &\leq& c_2 + c_3 {\rm Re} \, \zeta_m \leq c_4 + c_5 \log^+ |X_m| + c_6 \log |\zeta_m| ,\\
|\zeta_m| &\leq& c_7 + c_8 \log^+ |X_m|  \leq c_9 + c_{10} \log m. 
\end{eqnarray*}
Now (\ref{p4}) is obtained by applying the Koebe distortion theorem \cite[Theorem 1.3]{Hay9} with  
$$\mu ( \lambda ) = \phi(v), \quad v =  \frac{1+\lambda}{1-\lambda} , \quad | \lambda | < 1, \quad \zeta_m = \frac{1+\lambda_m}{1-\lambda_m} .
$$
This yields, for large $m$, since $\zeta_m$ tends to infinity in $K_3$, 
$$
| \phi(\zeta_m)| = | \mu (\lambda_m) | \leq \frac{c_{11} }{\left(1 - | \lambda_m|^2 \right)^2} =
c_{11} \left| \frac{| \zeta_m|^2 + 2 {\rm Re}\,  \zeta_m  + 1 }{4 {\rm Re}\,  \zeta_m  } \right|^2 \leq c_{12} | \zeta_m|^2 \leq c_{13} ( \log m )^2 . 
$$
%integrating $\phi'/\phi$ along the  line segment from $v_0 $ to $\zeta_m^* = {\rm Re} \,  \zeta_m + i {\rm Im} \, v_0 $, followed by that from
%$ \zeta_m^*$ to $\zeta_m \in K_3$;  this yields, for large $m$, in view of (\ref{h3}),
%$$ | \phi(\zeta_m)| \leq c_{11} +  c_{12} |\zeta_m|^{4 \pi } \leq c_{13} + c_{14} (\log m)^{4 \pi} .%\leq  (\log m)^{5 \pi} .$$
\hfill$\Box$
\vspace{.1in}

\section{Proof of Theorem \ref{thmel}}\label{pfthme1}

Let $A$, $f_1 $ and $ f_2$  be as in the hypotheses, without loss of generality satisfying $W(f_1, f_2) = \pm 1$, and set $E = f_1f_2$.
%Then the zeros of $E$ have finite exponent of convergence $\lambda(E)$. 
Choose $M > 0$ such that $|A(0)| $ and all finite critical and asymptotic values of $A$ have modulus at most $M/2$. It may be assumed that
$M \leq 1$, because otherwise the $f_j$ may be replaced by the functions $g_j(z) = M^{1/2} f_j (z/M)$, which solve
$$
y'' + B(z) y = 0, \quad B(z) = M^{-2} A(z/M  ).
$$
%As shown in \cite{BIL1}, the assumption that $\lambda(f_j) < + \infty$ for $j=1,2$ and
%the Bank-Laine equation (\ref{bleq}) together imply that $E$ has finite order.  
If $T(r, E) = O(r)$ as $r \to + \infty$, then (\ref{bleq}) 
delivers (\ref{AEest}),
%$T(r, A) = O(r)$ as $r \to + \infty$, 
and Theorems \ref{thmA} and \ref{thmB} force $\rho(A) = \rho(E) = 1$. 

Assume henceforth that $T(r, E) \neq O(r)$ as $r \to + \infty$ and, following standard notation of the Wiman-Valiron theory \cite{Hay5},
denote by $\mu(r, E)$  the maximum term of the Maclaurin 
series of $E$, and by $\nu(r, E)$ the central index. 
Then inequalities from \cite{Hay5} give
\begin{equation}
 \label{hay5ineq}
T(r, E) \leq \log^+ M(r, E) \leq \log^+ \mu(2r, E) + \log 2 \leq  \int_1^{2r} \frac{\nu(t, E)}{t} \, dt + O(1) ,
\end{equation}
and so it may be assumed that $\nu(r) = \nu(r, E) \neq O(r)$ as $r \to + \infty$. 

Let $1/2 < \tau < 1$. 
It follows from the Wiman-Valiron theory \cite{Hay5} that there exists a sequence $(z_n)$ satisfying 
%$|z_n| = r_n \to + \infty$ and $|E(z_n)| = M(r_n, E)$ with
\begin{equation}
 \label{b0}
|z_n| = r_n \to + \infty, \quad |E(z_n)| = M(r_n, E), \quad 
 \lim_{n \to + \infty} \frac{\nu(r_n)}{r_n} = + \infty ,
\end{equation}
such that, if $z = z_n e^\sigma$, $|\sigma| < \nu(r_n)^{-\tau} $, then
%\begin{equation} \label{b1}
$$
 E(z) \sim \left( \frac{z}{z_n} \right)^{\nu(r_n)} E(z_n), \quad 
 \frac{E'(z)}{E(z)} \sim \frac{\nu(r_n)}{z} , \quad 
 \frac{E''(z)}{E(z)} \sim \frac{\nu(r_n)^2}{z^2} ,
$$
%\end{equation}
%for $z = z_n e^\sigma$, $|\sigma| < \nu(r_n)^{-\tau} $. Combining (\ref{bleq}) with (\ref{b1}) shows that 
as well as, in view of (\ref{bleq}), 
%\begin{equation} \label{b2}
$$
 A(z) \sim - \frac{\nu(r_n)^2}{4 z^2} , \quad  A(z)^{-1/2} \sim \pm   \frac{2iz}{\nu(r_n)} .
$$
Thus (\ref{b0}) delivers   $\min \{ |E(z_n)|, |A(z_n)| \} \to + \infty $  as $n \to + \infty$, while 
applying Cauchy's estimate 
for derivatives to $A^{-1/2}$ yields 
$$
 A(z_n)^{-3/2} A'(z_n) = O \left(  \frac{r_n}{\nu(r_n)} \cdot \frac{\nu(r_n)^\tau}{r_n} \right) \to 0.
$$
Take $n$ so large that
\begin{equation}
 \label{b3}
|E(z_n)| > 2, \quad 
 \log |A(z_n)| > N_1, \quad \left| \frac{A(z_n)^{3/2}}{A'(z_n)} \right| > N_1, 
\end{equation}
where $N_1$ is the positive constant from  Proposition \ref{mainprop}.
Then $z_n$ lies in a component $C$ of
$\{ z \in \C : \, |A(z)| > 1 \}$, and $0 \not \in C$ since $|A(0)| < 1$. 
Because 
all finite critical and asymptotic values of $A$ have modulus at most $1/2$,
a change of variables (\ref{b4}) 
gives a conformal equivalence between $C$ and the right half-plane ${\rm Re} \, v > 0$. 
Choose $\sigma_n$, with ${\rm Re} \, \sigma_n > 0$, such that 
$z_n = \phi(\sigma_n)$. 
Then $e^{\sigma_n} = A(z_n)$ and (\ref{b4}) and
(\ref{b3}) imply that
%, as $n \to + \infty$, 
\begin{equation}
 \label{b5} 
 {\rm Re} \, \sigma_n > N_1 , \quad \left| e^{\sigma_n/2} \phi'(\sigma_n) \right| = \left| \frac{A(z_n)^{3/2}}{A'(z_n)} \right| 
 > N_1 . 
\end{equation}
A solution $y(z)$ of (\ref{de1}) transforms 
under  (\ref{b4}) to a solution $w(v) = y(z)$ of (\ref{p2}), and $\{ f_1, f_2 \}$ to a pair of solutions
 $\{ w_1, w_2 \}$ of (\ref{p2}) with $W(w_1, w_2) = \pm \phi'$. Use (\ref{b5}) to apply Proposition \ref{mainprop}, with $v_0 = \sigma_n = \phi^{-1} (z_n)$.
Since $|w_1(v_0)w_2(v_0)| = |E(z_n)| > 2$, by (\ref{b3}), the function $E$ has a sequence of distinct zeros
$\phi(\zeta_m)$  satisfying (\ref{p4}). But this gives $c > 0$ such that, for all large $m \in \N$, and for 
$r$ satisfying $c (\log m)^{2 } \leq r < c (\log (m+1))^{2 }$,
$$
 n(  c (\log m)^{2 } , 1/E) \geq \frac{m}2  , \quad n(r, 1/E) \geq \frac{m+1}3 > \frac13 \exp \left( (r/c)^{1/2} \right), 
$$
which establishes (\ref{Ezerosconc}) and completes the proof of Theorem \ref{thmel}. 
\hfill$\Box$
\vspace{.1in}

The following example shows that the exponent $1/2$ in (\ref{Ezerosconc}) is sharp. Let $A(z) = \cos \sqrt{z}$, which belongs to the Eremenko-Lyubich class
$\mathcal{B}$, and let $f$ be a non-trivial solution of (\ref{de1}). Let $\nu(r, f)$ be the central index of $f$ 
and apply to $f$ the same results from the Wiman-Valiron theory \cite{Hay5} as used in
(\ref{hay5ineq}) and subsequently. If $r$ is large and lies outside an exceptional set of finite logarithmic measure, and if $|z| = r$ and $|f(z)| = M(r, f)$,
then 
$$
\frac{\nu(r,f)^2}{z^2} \sim \frac{f''(z)}{f(z)} = - A(z), \quad \nu(r, f) \leq \exp( c \sqrt{r} ),
$$
for some positive constant $c$. This upper bound for the non-decreasing function
$\nu(r, f)$ then holds for all large $r$, possibly with a larger $c$, and so applying to $f$
the inequalities  of (\ref{hay5ineq}) gives
$d_1 > 0 $ with 
$$
n(r, 1/f) \leq N(3r, 1/f) \leq T(3r, f) + O(1) \leq \exp ( d_1 \sqrt{r} )  
$$
as $r \to + \infty$. Because $\rho(A) = 1/2$, conclusion (A) of Theorem \ref{thmel} cannot hold in this case, and so the exponent $1/2$ in 
(\ref{Ezerosconc}) is sharp. 
\hfill$\Box$
\vspace{.1in}

% in that 
\section{Proof of Theorem \ref{thmBLpositive}}

Let $E$ be as in the hypotheses, and assume that the zeros of $E$ have exponent of convergence
$\lambda \leq 3/2$. Then   the canonical product $\Pi_0$ over these zeros has order  $\lambda$, and  
$E = \Pi_0 \exp( P_0 )$, with $P_0$ a real polynomial. 
%Then  $P_0 + iP_1$ has degree at most $1$, and so has $P_1$. 
%Since $\exp(iP_1) = \pm 1$ at every zero of $E$, it may be assumed that $P_1 = 0$ and $E$ is real entire. 
If $P_0$ has degree at least $2$, then 
the zeros of $E$ have Nevanlinna deficiency $\delta(0, E) = 1$, which contradicts \cite[Theorem 4.1]{qcsurg} 
(see also \cite[Theorem 2.1]{Lasparse}). It may therefore be assumed that $E$ 
has order  $\lambda \leq 3/2$.

There exist an entire function $A$ and solutions $f_1, f_2$ of (\ref{de1}) such that 
$W(f_1, f_2) = 1$ and $ E = f_1 f_2$.  Then $f_j(z) = 0$ gives $E'(z) = (-1)^j$ and
each $f_j$ has infinitely many zeros, as may be seen by  considering the graph of $E$ on the real axis. 
Define $U$ by
\begin{equation*}
 \label{Udef1a}
U = \frac{f_2}{f_1}, \quad \frac{U'}{U} = \frac{W(f_1, f_2)}{f_1 f_2} = \frac1E.
\end{equation*}
%Then $U$ is locally univalent and the fact that

%, and cannot be constant by elementary considerations. 

\begin{lem}
 \label{lemAtrans}
The  coefficient function $A$  in (\ref{de1}) has order at most $\lambda$ but is transcendental.
%The function $A$ .
\end{lem}
\textit{Proof.}
The first assertion is an immediate consequence of the Bank-Laine equation (\ref{bleq}).
The second may be deduced from a theorem of Steinmetz \cite{Steiradial}, or from a combination of Theorem~\ref{thmDL} with the result of
Edrei, Fuchs and Hellerstein \cite{EFL} that if $E$ is an entire function of finite order and genus at least $1$, all of whose zeros are
positive, then $0$ is a Nevanlinna deficient value of $E$, from which the transcendence of $A$ follows using (\ref{bleq}). 
It may also be proved using Hille's method as follows. Suppose that $A$ is a polynomial.
Since the $f_j$ have infinitely many positive zeros, the positive real axis must be one of the $2 + \deg (A)$  critical rays 
for the equation, and each $f_j$ must be large in both  adjacent sectors.  
%However, $A$ evidently cannot be constant, so l
Let $L$ be the first other critical ray encountered when
moving counter-clockwise from the positive real axis. Since the $f_j$ have only positive zeros, both $f_j$ must change from large to  small as this critical ray $L$ is crossed.
A contradiction then arises from the fact that linearly independent solutions cannot be small in the same sector, because the Wronskian is a non-zero constant.
\hfill$\Box$
\vspace{.1in}

Because $U'/U$ has order at most $3/2$ and is never $0$, while all zeros and poles of $U$ are simple, 
$U$ has no critical values and finitely many asymptotic values \cite{Lasing2016}.
%In particular, the inverse function $U^{-1}$ has finitely many singular values. 
Since $U'/U$ is real, there exists $\theta \in \R$ such that $U = f_2/f_1 = e^{ 2 i \theta } U_0$, with $U_0$ real meromorphic. But
replacing $f_1 $ by $f_1 e^{i \theta}$ and $f_2 $ by $f_2 e^{- i \theta}$ leaves $E$ unchanged, and so it may be assumed that $\theta=0$ and $U$ is real meromorphic.

Take  zeros $x_0 , x_1, x_2 \in \R$ of $f_2$, with $x_0 < x_1 < x_2$,
and let $R$ be the supremum of all $r > 0$ such that the branch of $U^{-1}$ mapping $0$ to $x_1$ admits unrestricted analytic
continuation in the open disc $B(0, r)$ of centre $0$ and radius $r$. Then $R$ is finite, and $U$ maps a simply connected domain $\Omega_1$, with $x_1 \in \Omega_1$, univalently onto
$B(0, R)$. Moreover, $U^{-1}$ has a singularity over some $\alpha$ with $|\alpha | = R$, and $\Omega_1$ contains a path $\gamma$ which  tends to infinity and
is mapped by $U$ onto the half-open line segment $[0, \alpha)$. If $\alpha $ is real then, because $U$ is real meromorphic and univalent on $\Omega_1$, 
the path  $\gamma $ must coincide with  $(-\infty, x_1]$ or $[x_1, + \infty)$, 
contradicting the fact that $x_0, x_2 \not \in \gamma$. 
Hence 
%The intersection of $\gamma$ with $\R$ is bounded, since $\Omega_1$ is simply connected and symmetric with respect to $\R$ but contains neither $x_0$ nor $x_2$. 
%Thus the symmetry of $\Omega_1$ and the univalence of $U$ on $\Omega_1$ together imply that 
$\alpha \not \in \R$ and, since $U$ has finitely many critical and asymptotic
values,  $U^{-1}$ has logarithmic singularities over $\alpha$ and $\overline{\alpha}$.
% but no other singularities over finite non-zero  values. 

\begin{lem}
 \label{lemsingE}
Let $F(z) = (E(z)-E(0))/z$. Then there exist $M_0 > 0$ and disjoint non-empty components $\Sigma_1, \Sigma_2$ of the set $\{ z \in \C : \, |F(z)| > M_0 \}$.
%The inverse of the function  has at least two direct singularities over $\infty$.
\end{lem}
\textit{Proof.} There exists $M > 0$ such that  for each  $\beta \in \{ \alpha, \overline{\alpha} \} $ there is   a component $\Omega_\beta $ of the set 
$\{ z \in \C : \, | U(z) - \beta | < 1/M \}$ mapped univalently by $v = \log 1/(U(z)-\beta)$  onto the half-plane $H_0$ given by ${\rm Re} \, v > \log M$. 
It may be assumed that $M$ is so large that $\Omega_\beta \cap B(0, 1) = \emptyset$ and  $\Omega_\alpha \cap \Omega_{\overline{\alpha}} = \emptyset$. 
Let $\phi: H_0 \to \Omega_\beta$ be the inverse function and write 
\begin{equation}
 \label{vdef1}
U(z) = \beta + e^{-v} , \quad z = \phi(v) \in \Omega_\beta , \quad v \in H_0 .
\end{equation}
Then 
\begin{equation}
 \label{Eform}
E(z) = \frac{U(z)}{U'(z)} = \frac{\beta  + e^{-v}}{-e^{-v}} \cdot \phi'(v)  = - (1+\beta e^v) \phi'(v),
\end{equation}
and $\phi$ satisfies, on $H_0$,  as in (\ref{h3}), 
\begin{equation}
 \label{h3a}
\left| \frac{\phi''(v)}{\phi'(v)} \right| \leq \frac4{{\rm Re} \, v - \log M}, \quad  \left| \frac{\phi'(v)}{\phi(v)} \right| \leq \frac{4 \pi}{{\rm Re} \, v - \log M}  .
\end{equation}
It follows from (\ref{vdef1}),
(\ref{Eform}) and (\ref{h3a}) that there exists $c_1 > 0$ such that, as $v \to + \infty$ on $\R$,
$$
|z| = |\phi(v)| = o( v^{c_1} ) = o( e^v | \phi'(v) | ) = o( |E(z)| ) , \quad F(z) \to \infty, 
$$
%and so $F(z)$ is unbounded on $\Omega_\beta$, 
whereas if ${\rm Re} \, v = 1 + \log M $ then 
$$
|E(z)| \leq |z| ( 1 + | \beta | M e )  \left| \frac{\phi'(v)}{\phi(v)} \right| \leq |z| ( 1 + | \beta | M e ) 4 \pi , \quad 
|F(z)| \leq (1 + | \beta | M e ) 4 \pi + |E(0)| .
$$
Hence there exist $M_0 > 0 $ such that the set 
$\{ v \in H_0: \, |F(\phi(v))| > M_0 \}$ has a component whose closure with respect to the finite plane lies in $H_0$.
%$\Omega_0 \subseteq \Omega_\beta  $ of the set $\{ z \in \C : \, |F(z)| > M_0 \}$.
% such that $\Omega_0 $ and its boundary lie in $\Omega$. 
\hfill$\Box$
\vspace{.1in}

The remainder of the proof follows lines fairly similar to \cite{Ros,Shen}. 
By (\ref{bleq}) and well known estimates for logarithmic derivatives \cite{Gun2}, there exist positive integers $M_1, M_2 $ 
such that 
\begin{equation}
 \label{Abound}
\left| \frac{E'(z)}{E(z)} \right| + \left| \frac{E''(z)}{E(z)} \right|
\leq |z|^{M_1} , \quad 
|A(z)| = \frac1{4 |E(z)|^2}  + O \left(  |z|^{M_2}  \right),
\end{equation}
for all $z$ outside a union $U_1$ of countably many open discs, whose centres tend to infinity and whose radii have finite sum. 
Choose a polynomial $P$, of degree at most $M_2$, such that
\begin{equation}
 \label{Bdef}
B(z) = \frac{A(z)-P(z)}{z^{M_2+1}} 
\end{equation}
is entire. 
For $j= 1, 2$ 
define a subharmonic function $u_j(z)$ on $\C$ by $u_j(z) = \log |F(z)/M_0|$ on $\Sigma_j$, with $u_j(z) = 0$ on $\C \setminus \Sigma_j$,
where $M_0$ and $\Sigma_1, \Sigma_2$ are as in Lemma \ref{lemsingE}. 
Similarly, let $\Sigma_3$ be a component of the set $\{ z \in \C : \, |B(z)| > 1 \}$, and set
$u_3(z) = \log |B(z)|$ on $\Sigma_3$, with $u_3(z) = 0$ on $\C \setminus \Sigma_3$. These $u_j$ have orders satisfying
$\rho( u_j) \leq \rho (F) = \rho (E) = \lambda$, for $j=1, 2$, while $\rho(u_3) \leq \rho(B) = \rho(A) \leq \lambda $. 

For $j = 1, 2, 3$ and $t > 0$ let $\theta_j(t)$ be the angular measure 
of the intersection of $\Sigma_j $ with   the circle $S(0, t)$ of centre $0$ and radius $t$.
If $j = 1, 2$ and $|z|$ is large and $z \in \Sigma_j \cap \Sigma_3$, then (\ref{Bdef}) implies that
$z$ lies in the exceptional set $U_1$ of (\ref{Abound}). Hence there exists a set 
$F_0 \subseteq [1, + \infty )$, of
finite linear measure, such that if $r \in [1, + \infty) \setminus F_0$  then the following all hold: (a)   $S(0, r)$
does not meet $U_1$; (b)  
$\Sigma_j \cap \Sigma_{j'} \cap S(0, r) = \emptyset $ for $j \neq j'$; (c) no $\Sigma_j$ contains  $S(0, r)$. 

Let $S$ be large and positive: then
%,  since each $u_j $ has order less than $\lambda < 3/2$, 
a well known consequence of Carleman's 
estimate for harmonic measure \cite[pp.116-7]{Tsuji} gives, as $r \to + \infty$,
\begin{eqnarray*}
 9 \log \frac{r}{S} &\leq & \int_{[S, r] \setminus F_0}  \left( \sum_{j=1}^3 1 \right)^2 \, \frac{dt}t + O(1)  \leq
 \int_{[S,r] \setminus F_0}  \left( \sum_{j=1}^3 \theta_j(t) \right) \left( \sum_{j=1}^3 \frac1{\theta_j(t)} \right) \, \frac{dt}t + O(1) \\
&\leq& 2 \sum_{j=1}^3 \int_{[S,r] \setminus F_0}   \frac\pi{t \theta_j(t)}  \, dt  + O(1)  \leq 2 \sum_{j=1}^3 \log ( \max \{ u_j(z): \, |z| = 2r \}  )  + O(1) 
\\
&\leq& 2   \sum_{j=1}^3 (\rho (u_j) + o(1)) \log r \leq ( 6 \lambda + o(1)) \log r \leq ( 9 + o(1)) \log r .
\end{eqnarray*}
It follows at once that 
%$$\frac92 \leq \sum_{j=1}^3 \rho (u_j) \leq 3 \lambda \leq \frac92,$$and so  
$\rho (u_j) = \lambda = 3/2$ for each $j$. 
\hfill$\Box$
\vspace{.1in}

%\section{A quasiconformal mapping}\label{example1}

\section{A Bank-Laine function  with positive zeros}\label{example3}

The construction of  an example demonstrating the sharpness of Theorem \ref{thmBLpositive} will involve 
domains $D_0$, $D_1$, $D_2$ and $D_3$ defined by  
\begin{eqnarray}
 \label{D1def}
D_0 &=& \{ u \in \C : \, 0 < |u| < + \infty , \, 0 < \arg u < 3 \pi /2 \} , \nonumber \\
D_1 &=&  
%D_0 \setminus  \{  s + it: \,  \pi /2 \leq s < + \infty , \, 0 < t < + \infty   \}, \nonumber \\
E_1 \cup E_2 ,\nonumber \\
E_1 &=&  \{  s + it: \, - \pi /2 < s < 0, \, - \infty < t < + \infty  \} , \nonumber \\
E_2 &=& \{  s + it: \, - \pi /2 < s < \pi/2 , \, 0 < t < + \infty  \} , \nonumber \\
D_2 &=& \{ v \in \C : \, 0 < |v| < + \infty, \, - \pi/2 < \arg v < 0 \} ,\nonumber \\
D_3 &=& D_2 \cup \{ \zeta \in \C : \, | \zeta | < 1, \, {\rm Re} \, \zeta > 0 \}.
\end{eqnarray}

%In the following lemma, the statement that a real-valued function $g$ has positive upper and lower bounds on an interval $I$ means simply that there exists 
%$\varepsilon > 0$ such that $\varepsilon < g(x) < 1/\varepsilon$ for all $x \in I$. 

\begin{lem}
 \label{lemqc2} 
Let $h: (-\infty, 1] \to (-\infty, 0] $ be a continuous bijection, such that $h(1)=0$ while $h'$ is continuous and has positive upper and lower bounds
for $- \infty < y < 1$ (that is, there exists 
$\varepsilon > 0$ such that $\varepsilon < h'(y) < 1/\varepsilon$ for $- \infty < y < 1$). Then there exists a homeomorphism $\psi $ from the closure of $D_3$ to that of 
$D_2$, such that: (A) $\psi$ maps $D_3$ quasiconformally onto $D_2$, with $\psi(z) \to \infty$ and $\psi(z) = O( |z| )$ as $z \to \infty$ in $D_3$; (B)
$\psi(iy) = i h(y) $ for $-\infty < y \leq 1$; (C) $\psi(z)$ is real and strictly increasing as $z  $ describes  
the boundary of $D_3$ clockwise from  $i$ to infinity. 
\end{lem}
\textit{Proof.} Let $\phi : D_3 \to D_2$ be a conformal bijection such that $\phi(i) = 0$ and $\phi(z) \to \infty$ as $z \to \infty$ in $D_3$. 
Then $\phi(z)$ is real and strictly increasing 
as $z  $ describes the boundary of $D_3$ clockwise from  $i$ to infinity.  Moreover, there exists
a continuous bijection $k: (-\infty, 1] \to (-\infty, 0] $ with $k(1) = 0$ and $\phi(iy) = i k(y)$. It is clear from the reflection principle that 
$k'(y) $ is continuous and positive  for $-\infty < y < 1$, and it will be shown  that $k'(y)$ has positive upper and lower bounds for $-\infty < y < 1$.

Take the restriction of $\phi$ to $\{ z \in D_3: \, |z| > r_1 \}$, for some large positive $r_1$, 
and reflect twice, first across the imaginary axis and then across the real axis. This shows that 
$$L_1 = \lim_{y \to - \infty } k'(y) = \lim_{z \to \infty} \phi'(z) $$ 
exists and is finite and positive, and gives 
$\phi(z) = O( |z| )$ as $z \to \infty $ in $D_3$. 

Next, extend $\phi : D_3 \to D_2$ by reflection across the imaginary axis 
to a conformal mapping onto the lower half-plane, and apply the reflection principle to $\phi_1(u) = \phi(e^{iu} )$ 
on the half-disc $\{ u \in \C : \, | u- \pi/2 | < r_2, \, {\rm Im } \, u > 0 \}$, for some small positive $r_2$. This extended function has $\phi_1'(\pi/2) \neq 0$, 
which shows that $L_2   = \lim_{y \to 1-} k'(y) $ exists and is finite and non-zero, and hence positive by continuity. 

The function $H = h \circ k^{-1}$ is a continuous bijection from $(-\infty, 0]$ to itself, and so there exists
a homeomorphism  $\eta$ from the closure of $D_2$ to itself given by 
$\eta( x+iy) = x + iH(y)$ for $x \geq 0$ and $y \leq 0$.  
Furthermore, the chain rule shows that $H'(y)$ is continuous, with positive upper and lower bounds, for $-\infty < y < 0$.
Hence $\eta$ is $C^1$ on $D_2$ with  
$$
2 \frac{\partial \eta}{\partial \overline{z}} = \eta_x + i \eta_y =  1 - H'(y)  , \quad  2 \frac{\partial \eta}{\partial z} = \eta_x - i \eta_y = 1 + H'(y)  ,
$$
which ensures that $\eta$ is quasiconformal on $D_2$, 
and that $\eta(z) = O( |z| )$ as $z \to \infty $ in $D_2$. 

It now follows that $\psi = \eta \circ \phi$ is a homeomorphism  from the closure of $D_3$ to that of $D_2$,
quasiconformal on $D_3$ itself, and satisfies
$$
\psi(iy) = \eta ( \phi(iy)) = \eta ( i k(y)) = i H( k(y)) = i h(y)  \quad \text {for $-\infty < y \leq 1$.}
$$
Finally, $\psi(z) = O( | \phi(z)| ) = O( |z|)$ as $z \to \infty $ in $D_3$. 
\hfill$\Box$
\vspace{.1in}

%\begin{figure}\includegraphics[width=16cm]{domains2.eps}\caption{The domains $D_1$ and $D_3$; $A^*, B^*, C^*, D^*$ are the images of $A, B, C, D$ under $u \rightarrow e^{iu}$}
%\label{figdomains}\end{figure}
%\section{A locally univalent function on a sector}\label{example2}

\begin{lem}
 \label{lemqrmap}
Let $E_0$ be the closure  of the domain $D_0$ in (\ref{D1def}), and define $F$ on $ E_0 \setminus D_1 $ by 
\begin{eqnarray}
 \label{Fprop1}
F(s +it) &=& f_1(s +it) \quad \text{for $- \infty < s \leq -\pi/2$,  $t \in \R$}, \nonumber \\
F(s +it) &=& f_2(s +it) \quad \text{for $\pi/2 \leq s < + \infty $,  $0 \leq t < + \infty $}, \nonumber\\
f_1(u) &=& - i \exp( 2 e^{iu} ), \nonumber\\
f_2(u) &=&  \cot (u/2) = -i \left( \frac{1 + e^{iu} }{1 - e^{iu}} \right) .
\end{eqnarray}
Then $F$ extends to a mapping from $E_0$ into the extended plane, continuous with respect to the spherical metric,
with the following properties.\\
(i) $H = \log F$ maps $D_1$ quasiconformally onto the quadrant $D_2$, with $H(\pi/2) = 0$.\\
(ii) Let $L_0$ be the path consisting of the line segment from $\pi$ to $0$ followed by the negative imaginary axis in the direction of
$- i \infty$. Then $F(u)$ is real and strictly increasing as $u$ describes $L_0$, and each $u_0 \in L_0$ has $s_0 > 0$ such that 
${\rm Im} \, F(u) < 0$ on $D_0 \cap B(u_0, s_0)$. \\
%, with $\log F(\pi/2) = 0$;\\
(iii) $F$ is locally injective on $E_0$; \\
(iv) There exists $c > 0$ such that $|F(u)| \leq \exp \exp ( c |u| )$ for $u \in D_0$ lying on the circles $|u| = (4n+1)\pi/2$, $n  \in \N$. 
\end{lem}
\textit{Proof.} 
Using the principal argument and logarithm, set 
\begin{eqnarray}
 h(y) &=& - \, \frac{\pi}2  + 2y  \quad \text{for $- \infty < y \leq 0$}, \nonumber \\
h(y) &=& - \, \frac{\pi}2  + \arg \left( \frac{1 + iy }{1 -iy  } \right) = - \, \frac{\pi}2  -i
 \log \left( \frac{1 + iy }{1 -iy  } \right) \quad \text{for $0 < y \leq 1$}.
\label{hmatch}
\end{eqnarray}
%with the argument chosen so that $h(1) = 0$. 
%As $y$ decreases from $1$ towards $0$, the point $(1+iy)/(1-iy)$ travels clockwise around the unit circle from $i$
For $0 < y < 1$ this gives 
%to $1$; thus $h$ is  a continuous bijection from $(-\infty, 1]$ to $(-\infty, 0] $. Moreover, as $y \to 0+$, 
$$
h(y) = - \, \frac{\pi}2  + 2 \arctan y, \quad 
h'(y) = \frac2{1+y^2} \to 2 \quad \text{as $y \to 0+$,}
%\frac{1 + iy }{1 -iy  } = 1 + 2iy + \ldots  , \quad \arg \left( \frac{1 + iy }{1 -iy  } \right) = 2y + \ldots ,
$$
and so $h'(0) = 2$. 
Thus $h$ is  a continuous bijection from $(-\infty, 1]$ to $(-\infty, 0] $ and
%, by L'H\^opital's rule,  
$h'$ exists and is continuous  for $-\infty < y < 1$, with $1 \leq h'(y) \leq 2$ there. Lemma \ref{lemqc2} gives 
a homeomorphism $\psi $ from the closure of $D_3$ to that of $D_2$,  such that $\psi$ maps $D_3$ quasiconformally onto $D_2$, 
with $\psi(z) = O( |z| )$ as $z \to \infty$ in $D_3$ and $\psi(iy) = i h(y) $ for $-\infty < y \leq 1$. Furthermore, $\psi(z)$
is real and strictly increasing
as $z  $ describes the boundary of $D_3$ clockwise from  $i$ to infinity, and so is
%, by Lemma \ref{lemqc2}(C). 
$G = \exp \circ \psi$, which is continuous  on the closure of $D_3$ and satisfies, by (\ref{hmatch}),  
\begin{eqnarray}
G(v)  &=&  \exp ( i h(y) ) = - i \exp( 2 iy  ) = -i \exp( 2v ) \quad \text{for $v = iy$, $- \infty < y \leq 0$}, \nonumber \\
G(v)  &=& \exp ( i h(y) ) =  -i \left( \frac{1 + iy }{1 -iy  } \right) = -i \left( \frac{1 + v }{1 - v } \right) \quad \text{for $v=iy$, $0 < y \leq  1$}.
\label{Gmatch}
\end{eqnarray}

%the boundary of $D_3$ is described clockwise, starting from the point $i$. 

The next step is to set $F(u) = G(e^{iu}) $ on the closure of $D_1$. Now $v = e^{iu}$ 
maps $D_1$ conformally onto $D_3$, with  $v \to 0$ as
${\rm Im} \, u \to + \infty$ and $v \to \infty$ as $ {\rm Im} \, u \to - \infty$.
Furthermore, the boundary of $D_1$  is mapped by $v = e^{iu}$ as follows: the line ${\rm Re} \, u = - \pi/2$  to the negative imaginary axis;
the half-line ${\rm Re} \, u =  \pi/2$, $0 \leq {\rm Im} \, u < + \infty$,  to the segment 
$v=iy$, $0 < y \leq 1$; the  real interval $[0, \pi/2]$  to the arc of the unit circle from $1$ to $i$; the negative imaginary axis  to the real
interval $(1, + \infty)$.
Hence (\ref{Fprop1}) and (\ref{Gmatch}) imply that
%, coupled with the elementary properties of the mapping $e^{iu}$ from $D_1$ to $D_3$ illustrated in Figure \ref{figdomains}, imply that 
$F$ is well-defined and continuous
on  $E_0$, and that (i) and (ii) hold. Moreover, because $\psi$ is injective on $D_3$ and (\ref{Fprop1}) 
implies that each $\log |f_j(u)|$ is positive for $-\pi/2 < {\rm Re} \, u < \pi /2$ and negative for $\pi/2 < | {\rm Re} \, u | < 3 \pi /2$,
the function $F$ is locally injective  on $E_0$. Thus it remains only to prove (iv). By (\ref{Fprop1}), it is enough to bound the growth of $F(u)$ for $u \in D_1$, 
and hence it suffices to consider the continuous function $G(v) $ on the closure of $D_3$. But, 
as $v = e^{iu} \to \infty$ in  $D_3$,
$$
|F(u)| = |G(v)|\leq \exp( |\psi(v)| )  \leq \exp( O(|v| ) ) = \exp \left( O( |e^{iu} | )  \right) \leq \exp \exp ( 2 |u| ).
$$
\hfill$\Box$
\vspace{.1in}

Now define $V(z)$ on the open upper half-plane $H^+$  by $V(z) = F(z^{3/2})$, 
in which $z^{3/2}$ is the principal branch and 
$F$ is as in Lemma \ref{lemqrmap}. Then $V$ extends first to a (spherically) continuous function from the closed upper half-plane
into the extended plane,  mapping $\R$ into $\R \cup \{ \infty \}$, and then to the whole plane via
$V ( \overline z ) = \overline{V(z)}$. 
The extended function $V$ is locally injective on $\C$, in view of Lemma \ref{lemqrmap}, 
%locally injective on $\C \setminus I$, where $I = (-\infty, (\pi/2)^{2/3}]$. But if $x_0 \in I$ then Lemma \ref{lemqrmap} gives $s_0 > 0$ such that ${\rm Im} \,  V(z) < 0$ on
%$B(x_0, s_0) \cap H^+$ and $V(x)$ is strictly decreasing on $(x_0-s_0, x_0+s_0)$. 
%,  by Lemma \ref{lemqrmap}(i), (ii) and (\ref{fjprop}). Hence $V$ is , by  (\ref{fjprop}).
and quasimeromorphic, by \cite[Ch. I, Theorem 8.3]{LV}.
%or by noting  that the local integrability of partial derivatives  follows easily from  (\ref{h'est}), (\ref{gqcdef}) and the 
%reflection principle applied to $\psi^{-1}$, while $\psi^{-1}(z^{3/2}) $ has an analytic extension to neighbourhoods of
%$0$ and $ (\pi/2)^{2/3}$, by the reflection principle and (\ref{prop1a}) respectively. 
%: here %the whole plane by reflection: here one may appeal to Rickman p.179 but in this example the required absolute continuity on almost all horizontal and vertical lines 
%follows easily from (\ref{h'est}) and (\ref{gqcdef}). 
Lemma \ref{lemqrmap}(iv) delivers, as $n \to + \infty$,  
\begin{equation}
 \label{rndef}
\log^+ \log^+ |V(z)| =  O( |z|^{3/2} ) \quad \text{for} \quad |z| = r_n = \left( \frac{(4n+1)\pi}2 \right)^{2/3} . 
\end{equation}

The remainder of the construction proceeds much as in \cite{bebl1}. 
Let $D_4$ be the pre-image in $H^+$ of the domain $D_1$ under $u = z^{3/2}$, let $E_4$ be its closure, and $F_4$ the union of $E_4$ and its reflection across the real axis. 
Then $V$ is meromorphic off $F_4$ and writing $z =  re^{i \theta }$ and $u = se^{i \eta }$ 
shows that  the complex dilatation $\mu_V$ of $V$ satisfies, for some $C_1 , C_2 > 0$, 
\begin{eqnarray}
 \int_{1 \leq |z| < + \infty} \left| \frac{\mu_V(z)}{z^2} \right| \, dx dy  &\leq& 
%2 \int_{1 \leq |z| < + \infty, z \in D_4} \frac{1}{|z|^2} \, |dz|^2 \nonumber \\ &=&
2 \int_{1 \leq |z| < + \infty, z \in D_4} \frac{1}{r} \, dr \, d \theta  \nonumber \\
&=& C_1 \int_{1 \leq |u| < + \infty, u \in D_1} \frac{1}{s} \, d \eta \,  ds  \nonumber \\
&\leq& C_2\int_1^{+ \infty}  \frac{1}{s^2} \, ds = C_2 . 
\label{teich}
\end{eqnarray}
Let $\phi$ be the (unique) quasiconformal homeomorphism of the extended plane which solves the Beltrami equation 
$\displaystyle{ \frac{\partial \phi}{\partial \overline{z}} = \mu_V(z) \, \frac{\partial \phi}{\partial z} }$ a.e.
%$\phi_{\overline z} = \mu_V(z) \, \phi_z $ and  %with coefficient $\mu_V$, 
and fixes each of $0$, $1$ and $\infty$ \cite{LV}. In view of (\ref{teich}) and the Teichm\"uller-Belinskii theorem \cite[Ch. V, Theorem 6.1]{LV},
there exists $\alpha\in \C \setminus \{ 0 \}$ with 
\begin{equation}
 \label{phiasym}
\phi(z) \sim \alpha z
\end{equation}
as $z \to \infty$. Moreover, there exists a 
locally univalent meromorphic function $U$ such that $V = U \circ \phi$ on $\C$. Let $U_1(z) = \overline{U(\overline{z})}$: then 
$$
U(\phi(z)) = V(z) = \overline{V(\overline{z})} =  \overline{U(\phi(\overline{z}))} = U_1 \left( \overline{\phi(\overline{z})} \right) .
$$
Hence $\phi(z)$ and $\overline{\phi(\overline{z})}$ have the same complex dilation a.e. and both fix $0$, $1$ and $\infty$, which implies,
by uniqueness, that
$\phi$ is real on $\R$ and $U$ is real meromorphic. Furthermore, $\phi ([0, + \infty)) = [0, + \infty)$ and
all zeros and poles of $U$ are real and positive, while $E = U/U'$ is a real Bank-Laine
function with positive zeros.

Now $U$ satisfies, by Lemma \ref{lemqrmap} and (\ref{phiasym}), 
\begin{equation}
 \label{zeropoleest}
n(r, 1/U) + n(r, U) = O( r^{3/2} ) \quad \text{ as $r \to + \infty$. }
\end{equation}
Let $\Pi_1$ and $\Pi_2$ be the canonical products over the zeros and poles of $U$ respectively, which have order at most $3/2$, by
(\ref{zeropoleest}), and write
\begin{equation}
 \label{hdef}
U = \frac{\Pi_1}{\Pi_2} \, e^h, \quad 
\frac1E = \frac{U'}{U} = \frac{\Pi_1'}{\Pi_1} -  \frac{\Pi_2'}{\Pi_2} + h' ,
\end{equation}
where $h$ is an entire function. For $|z| = r_n$, where $n$ is large,  
(\ref{rndef}) and (\ref{phiasym}) yield 
$$\log^+ \log^+ |U (\phi(z))|= O( |z|^{3/2} ) = O( |\phi(z)|^{3/2} ).$$
%on the image under $\phi$ of the circles $|z| = r_n$, and t
Thus the maximum principle delivers 
$$
\log^+ \log^+ |\Pi_2(\zeta)U (\zeta)| = O( |\zeta|^{3/2} ) 
$$
as $\zeta \to \infty$ and hence 
\begin{equation}
 \label{PiUest}
\log T(r, \Pi_2 U) = O(r^{3/2}) \quad \text{as $r \to + \infty$}.
\end{equation}
Combining (\ref{PiUest}) with (\ref{hdef}) and the lemma of the logarithmic derivative leads to
$$m(r, h') \leq m \left(r, \frac{(\Pi_2 U)'}{\Pi_2U} \right) + m \left(r, \frac{\Pi_1'}{\Pi_1} \right) + O(1) 
= O( r^{3/2} ) \quad \text{as $r \to + \infty$}.
$$
Hence $h'$ and $E$ have order of growth at most $ 3/2$. Applying Theorem \ref{thmBLpositive} then shows that $E$ is a real Bank-Laine function 
whose zeros are all real and positive and have exponent of convergence $3/2$,
and that $E$ itself has order
$ 3/2$, as has the associated coefficient function $A$. 
%Thus $E$ provides a further counter-example to the Bank-Laine conjecture. 
\hfill$\Box$

\noindent
School of Mathematical Sciences, University of Nottingham, NG7 2RD.\\
james.langley@nottingham.ac.uk


\vspace{.1in}


{\footnotesize
\begin{thebibliography}{99}
\bibitem{BIL1}
S. Bank and I. Laine, On the oscillation theory of
$f''	  + Af = 0$ where $A$ is entire, Trans. Amer. Math. Soc. 273
(1982), 351-363.
%\bibitem{BIL2}S. Bank and I. Laine, Representations of solutions of periodic second order linear differential equations, J. reine angew. Math. 344 (1983), 1-21.
\bibitem{BIL3} S. Bank and I. Laine, On the zeros of meromorphic solutions of second-order linear differential equations, Comment. Math. Helv. 58 (1983), 656-677.
%\bibitem{BLL1} S. Bank, I. Laine and J.K. Langley, On the frequency of zeros of solutions of
%second order linear differential equations, Result. Math. 10 (1986), 8-24.
%\bibitem{BL1} S. Bank and J.K. Langley, On the oscillation of solutions of certain linear differential equations in the complex domain, Proc. Edinburgh Math. Soc.
%30 (1987), 455-469.
%\bibitem{BL2} S. Bank and J.K. Langley, Oscillation theory for higher order linear differential equations with entire coefficients, Complex Variables 16 (1991), 163-175.
\bibitem{Ber4}
W. Bergweiler, Iteration of meromorphic functions, Bull. Amer. Math. Soc. 29 (1993), 151-188.
%\bibitem{BCL} W. Bergweiler, J. Clunie and J.K. Langley, Proof of a conjecture of Baker concerning the distribution of fixpoints, Bull. London Math. Soc. 27 (1995), 148-154.
\bibitem{BE}
W. Bergweiler and A. Eremenko, On the singularities of the inverse to
a meromorphic function of finite order, Rev. Mat. Iberoamericana 11 (1995), 355-373.
\bibitem{bebl1}
W. Bergweiler and A. Eremenko, On the Bank-Laine conjecture,  J. Eur. Math. Soc. 19 (2017), 1899-1909.
\bibitem{bebl2}
W. Bergweiler and A. Eremenko, Quasiconformal surgery and linear differential equations, to appear, J. Analyse Math. 
\bibitem{qcsurg}
D. Drasin and J.K. Langley, Bank-Laine functions via quasiconformal surgery, Transcendental Dynamics and Complex Analysis, London 
Mathematical Society Lecture Notes 348 (2008), Cambridge University Press, 165-178. 
\bibitem{EFL}
A. Edrei, W.H.J. Fuchs and S. Hellerstein, 
Radial distribution and deficiencies of the values of a meromorphic function,
Pacific J. Math. 11 (1961) 135--151. 
\bibitem{EL}
A. E. Eremenko and M.Yu. Lyubich, Dynamical properties of some classes
of entire functions, Ann. Inst. Fourier Grenoble 42 (1992),
989-1020.
\bibitem{Gun2}
G. Gundersen, Estimates for the logarithmic derivative
of a meromorphic function, plus similar estimates, J. London Math. Soc.
(2) 37 (1988), 88-104.
%\bibitem{Hay1}W.K. Hayman, Picard values of meromorphic functions and their derivatives,Ann. of Math. 70 (1959), 9-42.
\bibitem{Hay2}
W.K. Hayman, Meromorphic functions, Oxford at the Clarendon Press, 1964.
\bibitem{Hay5}
W.K. Hayman, The local growth of power series: a survey of the
Wiman-Valiron method, Canad. Math. Bull. 17 (1974), 317-358.
%\bibitem{Hay7}
%W.K. Hayman, Subharmonic functions Vol. 2, Academic Press, London, 1989.
\bibitem{Hay9}
W.K. Hayman, Multivalent functions, 2nd edition, Cambridge Tracts in
Mathematics 110, Cambridge University Press, Cambridge 1994.
%\bibitem{HKLRT}
%J. Heittokangas, R. Korhonen, I. Laine, J. Rieppo and K. Tohge,
%Complex difference equations of Malmquist type, 
%Comput. Methods Funct. Theory 1 (2001), 27-39.
%\bibitem{HMR1}
%S. Hellerstein, J. Miles and J. Rossi, On the growth of solutions of $f
%'' + gf' + hf = 0$, Trans. Amer. Math. Soc. 324 (1991), 693-706.
%\bibitem{HMR2}
%S. Hellerstein, J. Miles and J. Rossi, On the growth of solutions of certain
%linear differential equations, Ann. Acad. Sci. Fenn. Ser. A. I. Math. 17
%(1992), 343-365.
\bibitem{Hil1} E. Hille, Lectures on ordinary differential equations, Addison-Wesley, Reading, Mass., 1969.
\bibitem{Hil2}
E. Hille, Ordinary differential equations in the complex domain, Wiley,
New York, 1976.
%\bibitem{Hin1} A. Hinkkanen, Reality of zeros of derivatives of meromorphic functions, Ann. Acad. Sci. Fenn. 22 (1997), 1-38.
%\bibitem{Hin2} A. Hinkkanen, Zeros of derivatives of strictly non-real meromorphic functions, Ann. Acad. Sci. Fenn. 22 (1997), 39-74.
\bibitem{Lai1}
I. Laine, Nevanlinna theory and complex differential equations,
de Gruyter Studies in Math. 15,  Walter de Gruyter, Berlin/New York 1993.
\bibitem{La5}J.K. Langley, Proof of a conjecture of Hayman concerning $f$ and $f'' $, J. London Math. Soc. (2) 48 (1993), 500-514.
%\bibitem{La9} J.K. Langley, On second order linear differential polynomials, Result. Math. 26 (1994), 51-82.
%\bibitem{Lades}J.K. Langley, On entire solutions of linear differential equations with one dominant coefficient, Analysis 15 (1995), 187-204.
%\bibitem{Laqc}J.K. Langley, Quasiconformal modifications and Bank-Laine functions,Archiv der Math. 71 (1998), 233-239.
%\bibitem{Laperm} J.K. Langley, Permutable entire functions and Baker domains, Math. Proc. Camb. Phil. Soc. 125 (1999), 199-202.
\bibitem{Lasparse}
J.K. Langley, Bank-Laine functions with sparse zeros, Proc. Amer. Math. Soc. 129 (2001), 1969-1978.
%\bibitem{Lades02} J.K. Langley, Linear differential equations with entire coefficients of small growth, Arch. Math. (Basel)78 (2002), 291-296.
%\bibitem{Larubel} J.K. Langley, Composite Bank-Laine functions and a question of Rubel, Trans. Amer. Math. Soc. 354 (2002), 1177-1191.
%\bibitem{Lanew}
%J.K. Langley,The second derivative of a meromorphic function of finite order, Bulletin London Math. Soc. 35 (2003), 97-108.
%\bibitem{Laoy} J.K. Langley, Integer points of entire functions,  Bulletin London Math. Soc. 38 (2006), 239-249. 
%\bibitem{schwarzian} J.K. Langley,  The Schwarzian derivative and the Wiman-Valiron property, to appear, J. Analyse Math.
\bibitem{Lasing2016}
J.K. Langley, Transcendental singularities for a meromorphic function with logarithmic derivative of finite lower order, to appear, 
Comput. Methods Funct. Theory. 
\bibitem{LV}
O. Lehto and K. Virtanen, Quasiconformal mappings in the plane,
2nd edn., Springer, Berlin, 1973.
\bibitem{Nev}
R. Nevanlinna, Eindeutige analytische Funktionen,
2. Auflage, Springer, Berlin, 1953.
\bibitem{Ros}
J. Rossi, Second order differential equations with transcendental coefficients,
Proc. Amer. Math. Soc. 97 (1986), 61-66.
\bibitem{Shen} L.C. Shen, Solution to a problem of S. Bank regarding the exponent of convergence of the solutions of a 
differential equation $f''	  + A f= 0$, Kexue Tongbao 30 (1985), 1581-1585.
\bibitem{Shen2}
L.C. Shen, Construction of a differential equation $y'' + Ay = 0$
with solutions having prescribed zeros, Proc. Amer. Math.
Soc. 95 (1985), 544-546.
\bibitem{sixsmithEL}
D.J. Sixsmith, Dynamics in the Eremenko-Lyubich class,
Conform. Geom. Dyn. 22 (2018), 185-224. 
\bibitem{Steiradial}
N. Steinmetz,
Linear differential equations with exceptional fundamental sets II,
Proc. Amer. Math. Soc. 117 (1993), no. 2, 355--358. 
\bibitem{Tsuji}
M. Tsuji, Potential theory in modern function theory,
Maruzen, Tokyo, 1959.

\end{thebibliography}
}
\end{document}